\documentclass[12pt,a4paper]{amsart}
\usepackage{graphicx}
\usepackage{epstopdf}
\usepackage{amssymb}
\usepackage{amsthm}
\usepackage{amsmath}
\usepackage{pstcol, pst-node}
\usepackage{mathrsfs}
\usepackage{amscd}

\usepackage{comment}
\usepackage[all]{xy}
\newtheorem{thm}{Theorem}

\newtheorem{cor}[thm]{Corollary}
\newtheorem{prop}[thm]{Proposition}
\newtheorem{lem}[thm]{Lemma}
\newtheorem{rem}[thm]{Remark}

\theoremstyle{definition}
\newtheorem{defn}[thm]{Definition}

\newtheorem{prop-def}[thm]{Proposition-Definition}

\newtheorem{example}[thm]{Example}

\newcommand{\Ext}{\mathop{\mathrm{Ext}}\nolimits}

\newcommand{\Hilb}{\mathop{\mathrm{Hilb}}\nolimits}

\newcommand{\res}{\mathop{\mathrm{res}}\nolimits}

\newcommand{\Spec}{\mathop{\mathrm{Spec}}\nolimits}

\theoremstyle{remark}

\begin{document}
\title{Obstructions to deforming maps from curves to surfaces}

\author{Takeo Nishinou}
\date{}
\thanks{email : nishinou@rikkyo.ac.jp}
\address{Department of Mathematics, Rikkyo University,
	3-34-1, Nishi-Ikebukuro, Toshima, Tokyo, Japan} 
\subjclass[2000]{}
\maketitle
\begin{abstract}
This paper studies the obstructions to deforming a map from a complex variety to another variety which is an immersion of codimension one. We extend the classical notion of semiregularity of subvarieties to maps between varieties, and show that it largely extends the applicability. Then we apply the main result to several situations. First, we give deformations of non-reduced curves on surfaces in a geometrically controlled way. Also, we give a simple but effective criterion for the vanishing of the obstructions to equisingular deformations of nodal curves on surfaces. Finally, we construct nodal curves with very small geometric genus on surfaces of general type. 
\end{abstract}

\section{Introduction}
In this paper, we study deformations of a map from 
 a complex variety to another variety which is an immersion of codimension one.
When the dimension of the target space is larger than two, we assume the domain is smooth. 
When the target space is a surface, 
 the domain curve is assumed to be reduced but otherwise 
 it can have any singularity.

The starting point of this study
 is the notion of \emph{semiregularity} which goes back to Severi \cite{S1, S2}.
Here a smooth complex curve $C$ on a smooth complex surface 
 $S$ is called semiregular if the canonical linear system of the 
 surface cuts out on $C$ a complete linear system (see Subsection \ref{subsec:sr}).
Severi proved that when $C$ is semiregular, then its deformation on $S$ is unobstructed.

The notion of semiregularity and the associated result were generalized by Kodaira-Spencer
　\cite{KS} to smooth divisors on higher dimensional complex manifolds.
Later, Bloch \cite{B} extended the notion of semiregularity to local complete intersection subvarieties,
 and related it to the smoothness of the Hilbert scheme at the corresponding point
 as well as to variation of Hodge structures.
In particular, the semiregularity of a local complete intersection subvariety guarantees the 
 vanishing of the obstructions to the deformation.
More recently, these ideas have been generalized from multiple points of view, 
 see \cite{BF, I, IM} to name a few.

Although these results are striking, often it is not easy to check whether a given subvariety is semiregular or not, 
 and even if we start from semiregular subvarieties, standard construction such as taking coverings tends to 
 break the condition of semiregularity.
On the other hand, since the obstruction to deformations is in principle determined by the information 
 of a neighborhood of the subvariety, it would not be too optimistic to expect that we can extend the 
 notion of semiregularity from the original cohomological (in other words, global) condition to  
 a more local one.

In this paper, we show that this is in fact the case, and we extend the notion of semiregularity to 
 maps rather than subvarieties (such an extension was also considered in \cite{I} for maps with smooth domains, 
 from a different perspective).
\begin{thm}\label{thm:1}
	Let $\varphi\colon V\to X$ be a map from a smooth variety (when $\dim X>2$)
	or a reduced complete complex curve (when $\dim X = 2$) 
	to a smooth complex variety $X$, which is 
	an immersion whose image is a divisor in $X$.
	Assume that $\varphi$ is semiregular.
	Then the map $\varphi$ is unobstructed in the sense that any first order deformation
	can be extended to arbitrary higher order. 
\end{thm}
This largely extends the applicability of semiregularity, especially when $\dim X = 2$. 
The point is that in addition to making many curves semiregular which are not so in the classical sense, 
 it also makes it possible to control the properties of deformed curves.
Namely, while classical theorems simply assure the existence of deformations (and there is little control
 of the deformed curves),
  in Theorem \ref{thm:1}, we regard a curve on a surface as the image of a map, so that its deformations are restricted by 
 the possibilities of the deformations of the domain of the map.
Thus, taking the domain curve suitably, we can obtain various (e.g., equisingular)  deformations with desired properties. 

In practical situations, starting from a few embedded curves on a surface for which the classical semiregularity condition holds, 
 we construct new curves by putting these together in a simple combinatorial way.
Then the constructed curves, seen as the images of suitable maps, often satisfy the (extended) semiregularity again, 
 and we can deform such curves on the surface.

The most simple case is 
 when the target $X$ has the trivial canonical sheaf.
In this case, any immersion $\varphi$ from a reduced curve is 
 semiregular, so unobstructed.
We give examples how to apply Theorem \ref{thm:1} in practice in this situation.
In particular, we consider the deformation of non-reduced curves in a geometrically controlled way, 
 which may be difficult to deal with by other methods.
 
We also consider the case where the target $X$ is a more general surface.
In this case, not all maps $\varphi$ as above are semiregular, but we give a general criterion for 
 the semiregularity of maps whose image is a nodal curve.
Namely, let $C$ be a nodal curve and $\varphi\colon C\to X$ be an immersion.
Assume that the image $\varphi(C)$ is an embedded nodal curve and let $P\subset \varphi(C)$
 be the set of nodes such that the preimage $\varphi^{-1}(p)$, $p\in P$, consists of two points.
Then we have the following. 
\begin{thm}\label{thm:intro2}
Assume that $\varphi(C)$ is semiregular in the classical sense.
Then the map $\varphi$ is semiregular if and only if for each $p_i\in P$, there is a first order deformation of 
$\varphi(C)$ which smoothes $p_i$, but does not smooth the other nodes of $P$. 
\end{thm}
Geometrically, this is, in a sense, the opposite to the classical \emph{Cayley-Bacharach property},
 see for example \cite{BHPV}.
Namely, we have the following. 
See Corollary \ref{cor:geomCB} for the precise statement.
\begin{cor}
Using the same notation as above, assume $\varphi(C)$ is sufficiently ample.
Then the map $\varphi$ is semiregular
 if for each $p_i\in P$, there is an effective divisor algebraically equivalent to $\varphi(C)$
which avoids $p_i$ but passes through all points in $P\setminus\{p_i\}$. 
\end{cor}

Finally, using a version of  Theorem \ref{thm:intro2} specialized to reducible curves, we construct 
 curves with very large number of nodes on general surfaces.
When the Kodaira dimension of the surface is non-positive, it is not difficult to construct 
 such curves, using Theorem \ref{thm:1} above for example.
In particular, in the proof of the classical Severi problem on the irreducibility of Severi varieties by Harris \cite{Ha}, 
 the key step is to show the existence of deformations of a given nodal curve to those with larger nodes, 
 and eventually to an irreducible nodal rational curve.
For surfaces with positive Kodaira dimension, one cannot expect the existence of such deformations.
See \cite{CS} for results of this direction.
However, we prove that any sufficiently ample smooth curve can be deformed into 
 a nodal curve with very large number of nodes.
Namely, we prove the following. 

\begin{thm}\label{thm:intro3}
	Let $X$ be a smooth complex projective surface with an effective canonical class.
	Let $L$ be a very ample class.
	Then, there is a positive number $A$ which depends on $L$ such that 
	for any positive integer $n$, the numerical class of $nL$ contains
	an embedded irreducible nodal curve $C$ whose geometric genus is less than $An$.
\end{thm}
When $L$ is sufficiently large, we can take $A = g(C_1) + K_X.L$, where 
 $g(C_1)$ is the genus of a general curve in the class $L$ and $K_X$ is the canonical class of $X$.
Since a general curve of the class $nL$ has genus of order $\frac{L.L}{2}n^2$, 
 the curve $C$ contains very large number of nodes. \\

There are many directions to which Theorem \ref{thm:1} should be extended.
In view of Bloch's result, the method of this paper might well generalize to maps 
 which are immersions but the domain can be 
 general local complete intersection.
%We can cover the case of hypersurfaces where the domain is smooth
% by the techniques developed here (see Subsection \ref{subsec:higherdim}).
Also, again in view of \cite{B}, one can consider the case where the target $X$ deforms, see \cite{N2}.
By adding some more calculation,
 the argument in this paper can be applied
 even in some cases where the ambient space degenerates,  
 see \cite{NY}.
The proof of Theorem \ref{thm:1} in this paper partly depends on transcendental method.
It will also be important to pursue a purely algebraic proof and investigate the positive characteristic case, 
 see \cite{M}.

\section{Obstruction to  curves on surfaces}	
First we deal with the case when $\dim X = 2$.
The case where $\dim X>2$ (and the domain is smooth) is treated later in Subsection \ref{subsec:higherdim}.
We think varieties in the complex analytic category.
Let $X$ be a smooth surface, not necessarily compact.
Let $C$ be a connected compact curve without embedded points.
We assume all irreducible components of $C$ are reduced, but otherwise $C$ can have any 
 singularity.
Let 
\[
\varphi\colon C\to X
\]
  be a map 
 which is an immersion.
Namely, at each point of $C$, there is a neighborhood of it
 such that the restriction of $\varphi$ to that neighborhood is an isomorphism onto the image. 
We will study deformations of $\varphi$.
\begin{rem}
Actually, we do not need to assume the variety $X$ to be smooth.
To study the deformation theory of the map $\varphi$, information of a neighborhood of 
 the image will suffice.
Similarly, the curve $C$ need not be compact or without boundary.
See also Remark \ref{rem:locality}.
%This observation will be important even if we will be basically interested in 
% the case where $X$ is actually smooth and closed.

%However,
% in actual problems, it is important to know the obstructions to deform.
%For that purpose, it will be desirable to have good control of the space of canonical sections,
% see Theorem \ref{}. 
\end{rem}

Deformation of the map $\varphi$ is controlled by 
 the hypercohomology groups
\[
\Ext^i(\varphi^*\Omega^1_X\to \Omega^1_C, \mathcal O_C).
\]
In the case mentioned above, the map $\varphi^*\Omega^1_X\to \Omega^1_C$ is surjective and
 the kernel of it is a locally free sheaf known as the conormal sheaf $\nu_{\varphi}^{\vee}$ of the map $\varphi$.
In particular, in this case the deformation of $\varphi$ is controlled by the usual sheaf cohomology groups.

\begin{defn}
We write the locally free sheaf $\mathcal Hom(\nu_{\varphi}^{\vee}, \mathcal O_C)$
 by $\nu_{\varphi}$ and call it the normal sheaf of the map $\varphi$.
\end{defn}

By definition of $\nu_{\varphi}$, we have the following.

\begin{prop}
The set of first order deformations of the map $\varphi$ is isomorphic to 
 $H^0(C, \nu_{\varphi})$. 
The obstruction to the existence of (the second and higher order as well) lifts lives in $H^1(C, \nu_{\varphi})$. \qed
\end{prop}

%\begin{rem}\label{rem:lift}
%Note that we are considering deformations of $\varphi$ relative to the base $B$.
%If the family $\pi\colon \mathcal X\to B$ is not trivial,  
% then in general there may not be a first order deformation 
% of $\varphi$ relative to the base.
%If there is one, then other first order deformations are obtained by perturbing it by elements of
% $H^0(V, \nu_{\varphi})$. 
%Since there is no canonical choice of a first order lift, the correspondence between 
% the set of first order lifts and
% the set $H^0(C, \nu_{\varphi})$ is not canonical.
%However, if the fibration $\mathcal X\to B$ is trivial, 
% there is always a first order　canonical lift, 
% namely, the map $\varphi$ itself (over $\Bbb C[t]/t^2$),
% and the correspondence between 
% the set of first order lifts and
% the set $H^0(C, \nu_{\varphi})$ becomes canonical.
%\end{rem}

Now, let $\omega_C$ be the dualizing sheaf of $C$ defined by 
\[
\omega_C = \varphi^*K_X\otimes \nu_{\varphi}.
\]
Here $K_{X}$ is the canonical sheaf of $X$.
By the Serre duality, the group $H^1(C, \nu_{\varphi})$ is isomorphic to the dual of 
 $H^{0}(C, \nu_{\varphi}^{\vee}\otimes\omega_{C})$.
In particular, we have the following.
\begin{prop}\label{prop:obst}
The obstruction to deform $\varphi$ belongs to the dual of the cohomology group
 $H^{0}(C, \varphi^*K_X)$.\qed
\end{prop}

\subsection{Relation to semiregularity}\label{subsec:sr} 
In the space $H^{0}(C, \varphi^*K_X)$, 
 there is a distinguished subspace $\mathcal V$ composed by those 
 coming  from  
 elements of $H^{0}(X, K_X)$.
When this subspace coincides with $H^{0}(C, \varphi^*K_X)$, the calculation of the obstruction
 becomes particularly simple.

In the case where $\varphi$ is an embedding, this condition is nothing but the 
 \emph{semiregularity} of the subvariety $\varphi(C)\subset X$ studied in \cite{S1, S2}
 and extended to higher dimensional cases in \cite{KS, B}.
In \cite{B}, the definition of semiregularity is extended to local complete intersection 
 subvarieties in any dimensional ambient space.
Namely, let $X$ be a smooth complex projective variety of dimension $n$, 
 and $Z\subset X$ be a local complete intersection of codimension $p$. 
Let $\mathcal I\subset\mathcal O_X$ be the defining ideal of $Z$.
Then the normal sheaf
\[
\mathcal N_{Z/X} = Hom_{\mathcal O_Z}(\mathcal I/\mathcal I^2, \mathcal O_Z)
\]
 is locally free.
Let $K_X = \Omega_X^n$ be the canonical sheaf, and define the invertible sheaves $K_{Z/X}$
 and $\omega_Z$ on $Z$ by
\[
K_{Z/X} = \wedge^p \mathcal N_{Z/X}^{\vee},\;\; \omega_Z = K_{Z/X}^{\vee}\otimes K_X.
\]

The natural inclusion $\varepsilon\colon \mathcal N^{\vee}_{Z/X}\to \Omega_X^1\otimes\mathcal O_Z$
 gives rise to an element
\[
\begin{array}{ll}
\wedge^{p-1}\varepsilon \in Hom_{\mathcal O_Z}\left(
   \wedge^{p-1}\mathcal N_{Z/X}^{\vee}, \Omega_X^{p-1}\otimes\mathcal O_Z \right)
    &= \Gamma\left((\Omega_X^{n-p+1})^{\vee}\otimes K_X\otimes K_{Z/X}^{\vee}\otimes\mathcal N_{Z/X}^{\vee}  \right)\\
   &= Hom_{\mathcal O_X}(\Omega_X^{n-p+1}, \omega_Z\otimes\mathcal N_{Z/X}^{\vee}).
\end{array}
\]

This induces a map on cohomology 
\[
\wedge^{p-1}\varepsilon\colon H^{n-p-1}(X, \Omega_X^{n-p+1})\to H^{n-p-1}(Z, \omega_Z\otimes \mathcal N_{Z/X}^{\vee}).
\]
In \cite{B}, the dual of this map 
\[
\pi\colon H^1(Z, \mathcal N_{Z/X})\to H^{p+1}(X, \Omega_X^{p-1})
\]
 is called the \emph{semiregularity map} and the subvariety $Z$ is called \emph{semiregular}
 if the semiregularity map is injective.

One of the main theorems of \cite{B} is the following, which generalizes results of \cite{KS, S1, S2}.
\begin{thm}\label{thm:B}
When $Z$ is semiregular in $X$, then the Hilbert scheme $\Hilb(X/\Bbb C)$ is smooth at the 
 point corresponding to $Z$.\qed
\end{thm}

In the case where $Z$ is a hypersurface so that $p = 1$, the above map $\wedge^{p-1}\varepsilon$
 on cohomology becomes
 the restriction $H^{n-2}(X, K_X)\to H^{n-2}(Z, i^*K_X)$, where $i\colon Z\to X$ is the inclusion.
In our situation where $n = 2$, while the map $\varphi\colon C\to X$ need not be an inclusion, the natural map
 $H^{0}(X, K_X)\to H^{0}(C, \varphi^*K_X)$ plays the similar role.
 
Although Theorem \ref{thm:B} is strikingly general, in explicit situations
 it is sometimes difficult to check the semiregularity, and divisors with higher multiplicities are
 usually not semiregular.  

In this paper, we remedy this point by local and direct calculation of obstructions.
In spite of the advances of the notions and technicalities concerning semiregularity
 in \cite{B} above and in even more sophisticated
 \cite{BF, I, IM}, the study of obstructions requires us to go back to 
 the classical method of \cite{KS}.

\section{Calculation of obstruction}\label{sec:sr}
As we mentioned in Subsection \ref{subsec:sr}, 
 when the natural map $H^{0}(X, K_X)\to H^{0}(C, \varphi^*K_X)$
 is surjective, then the map $\varphi$ will have a good property
 from the point of view of deformation theory.
\begin{defn}\label{def:semiregular}
We call the map $\varphi\colon C\to X$ semiregular
 if the natural map $H^{0}(X, K_X)\to H^{0}(C, \varphi^*K_X)$
 is surjective.
\end{defn}
However, to see this is a legitimate extension of the notion of classical semiregularity, 
 the cohomological method of previous works does not suffice, 
 and we need to develop a way to study obstruction more directly.

\subsection{Meromorphic differential forms and cohomological pairings}\label{subsec:paring}
We begin with giving a presentation of a \v{C}ech cocycle in a way suited to our purpose.
Let $C$ be a reduced curve and $\mathcal L$ be 
 an invertible sheaf on it.
Take an open covering $\{U_{1}, \dots, U_{m}\}$ of $C$ in the following way:
\begin{itemize}
	\item Each $U_{i}$ is a disk or an open subset which contains exactly one singular point of $C$.
	In the latter case, the normalization of $U_{i}$ is a union of disks where the number of the disks is the same
	as the number of the branches of $C$ at the singular point.
	\item For each singular point of $C$, there is a unique open subset from $\{U_{1}, \dots, U_{m}\}$
	containing it.
\end{itemize}
Here a disk means an analytic subset which is analytically isomorphic to the set 
 $D_r = \{z\in \Bbb C\;|\; |z|<r\}$ for some positive real number $r$.

If $U_i$ is the open subset containing a singular point $p_i$ of $C$, 
 let $U_i = U_{1,i}\cup \dots\cup U_{k_i,i}$ be the decomposition into branches.
The subsets $U_{i, j}$ and those $U_i$ which are disks form a covering of $C$ by locally closed subsets.
We write it by $\{V_j\}$.

\begin{rem}\label{rem:calldisk}
Although the normalization of $U_{i, j}$ is a disk, in general $U_{i, j}$ itself is not necessarily a disk.
However, we call it (and corresponding $V_j$) a disk for notational simplicity.
\end{rem}

By construction, we can naturally associate with each $V_j$ the unique
open subset from the covering $\{U_i\}$ above.
We write this open subset by $\tilde V_j$.
Namely, $\tilde V_j = U_i$ if $V_j = U_{a,i}$ for some $a$, 
and $\tilde V_j = U_i$ if $V_j = U_i$ is a disk.

We associate a local meromorphic section $\xi_j$ of $\mathcal L|_{V_j}$ with each $V_j$ in a way
 that if $V_i\cap V_j$ is an open subset, then $\xi_i-\xi_j$ is a holomorphic section of 
 $\mathcal L|_{V_i\cap V_j}$.

Let $V_1$ and $V_2$ be two of these locally closed subsets whose
intersection is an open subset (note that in this case $\tilde V_1\neq \tilde V_2$).
Then associate the section $\xi_1-\xi_2$ of $\mathcal L$ 
with the open subset 
$V_1\cap V_2$ (the order also matters, that is, we associate $\xi_2-\xi_1$
with $V_2\cap V_1$).
Taking each $U_i$ small enough, we can assume
\[
V_1\cap V_2 = \tilde V_1\cap \tilde V_2.
\]
Therefore, by associating $\xi_1-\xi_2$ with $\tilde V_1\cap \tilde V_2$, 
the set of sections $\{\xi_i\}$ determines a \v{C}ech 1-cocycle 
$\{\xi_{ij}\}$ with values in $\mathcal L$
for the covering $\{U_i\}$ above.

Conversely, any class in $H^1(C, \mathcal L)$ can be represented 
in this way (see Proposition \ref{prop:pair} below).\\

Let $\omega_{C}$ be the dualizing sheaf of $C$.
When $C$ is a reduced curve (as we are assuming in this paper), it is known that sections of $\omega_C$ are
 given by \emph{Rosenlicht differentials} (\cite{R}, see also \cite{AK, BHPV}).
Let $\nu\colon\widetilde C\to C$ be the normalization.
Let $\mathcal M_{\widetilde C}$ be the sheaf of meromorphic 1-forms on $\widetilde C$.
By definition, a section of $\nu_*\mathcal M_{\widetilde C}$ is called a
 meromorphic differential on $C$.
Then a meromorphic differential $\sigma$ is called a Rosenlicht differential
 if for all $x\in C$ and $g\in\mathcal O_{C, x}$, 
 \[
 \sum_{y_k\in\nu^{-1}(x)} \res(y_k, g\sigma) = 0
 \]
 holds.
Here $\res$ denotes the ordinary residue on $\widetilde C$.

%When $\{\eta_{ab}\}$ is seen as a \v{C}ech 1-cocycle of 
%$\mathcal L(\sum_{i=1}^k n_ip_i)$ with $n_i = -m_i$, 
Now let $\psi$ be a section of $\mathcal L^{\vee}\otimes \omega_{C}$.
Then $\{\xi_{ij}\}$ and $\psi$ make a natural pairing.
The value of this pairing is given as follows.
Namely, on a locally closed subset $V_i$, the fiberwise pairing 
 between $\xi_i$ and $\psi$ gives a meromorphic 
 section $\langle \psi, \xi_i\rangle$ of $\omega_C|_{V_i}$.
Let $p\in V_i$ be a point at which $\langle \psi, \xi_i\rangle$ has a pole, 
 and let $r_p$ be its residue.
If $p$ is contained in another $V_j$ and $V_i\cap V_j$ is an open subset
 (in this case, $p$ is not a singular point of $C$ and $\psi$ is a holomorphic 1-form there), 
 then since $\xi_i-\xi_j$ is holomorphic, 
 the residue of $\langle \psi, \xi_j\rangle$ at $p$ is the same as that of $\langle \psi, \xi_i\rangle$.
On the other hand, let $\{p_1, \dots, p_l\}$ be the set of singular points of $C$.
Let $V_{i_1}, \dots, V_{i_{a_i}}$
 be the disks on the branches at $p_i$ as above, and 
 let $r_{i_k}$ be the residue of $\langle \psi, \xi_{i_k}\rangle$ at $p_i$ on $V_{i_k}$
 (see Remark \ref{rem:calldisk} for our 
 convention of the terminology).
Note that, by definition, the residue of $\langle \psi, \xi_{i_k}\rangle$ at $p_i$ on $V_{i_k}$
 is the residue of the pull back of $\langle \psi, \xi_{i_k}\rangle$ to the normalization of 
 $V_{i_k}$ at the inverse image of $p_i$.
 
Then we have the following.
\begin{prop}\label{prop:pair}
	\begin{enumerate}
		\item Any cohomology class in $H^1(C, \mathcal L)$ can be represented
		by some set $\{\xi_i\}$ of local sections on locally closed subsets $\{V_i\}$ as above.
		\item 	
		The pairing between $\{\xi_{ij}\}$ and $\psi$
		is given by
		\[
		\langle \psi, \{\xi_{ij}\}\rangle = \sum_{i=1}^l\sum_{j=1}^{a_i} r_{i_j} + \sum_p r_p,
		\]
		where in the second summation, $p$ runs over the poles of $\langle \psi, \xi_i\rangle$
		for some $i$ which is not a singular point of $C$.
		This gives the natural nondegenerate pairing between $H^1(C, \mathcal L)$
		and its dual space $H^0(C, \mathcal L^{\vee}\otimes \omega_{C})$.
	\end{enumerate}
\end{prop}
\proof
We begin with proving that the given pairing is well-defined. 
This contains the following two claims.

First, we prove that the sum on the right hand side is independent of the choice of sections
$\{\xi_i\}$ which gives the same $\{\xi_{ij}\}$.
To see this, take another set of sections $\{\xi'_i\}$ which defines
the same \v{C}ech 1-cocycle.
Take an irreducible component $C_a$ of $C$.
Then, if two locally closed subsets $V_i$ and $V_j$ contained in $C_a$
 intersect and $V_i\cap V_j$ is open, the differences $\xi_i-\xi_i'$ and $\xi_j-\xi_j'$
coincide on $V_i\cap V_j$.
Thus, these differences define a global section $\xi''_a$ of 
 $\pi_a^*\mathcal L|_{C_a}$, where $\pi_a\colon \widetilde C_a\to C_a$ is
 the normalization.
Then the fiberwise pairing between $\xi''_a$ and the pull back of $\psi$ gives
a meromorphic 1-form on $\widetilde C_a$, and by the residue theorem
the sum of its residues is zero.
Doing this calculation on each component of $C$, 
 we see that the pairing $\langle \psi, \{\xi_{ij}\}\rangle$
 does not depend on the representative $\{\xi_i\}$.

Second, we need to show that the sum depends only on the cohomology class of 
$\{\xi_{ij}\}$. %, seen as a class of $H^1(\Sigma, \mathcal L(\sum_{i=1}^k n_ip_i))$.
Let $\{\xi_i\}$ and $\{\xi_i'\}$ be the sets of local sections defining the same 
cohomology class.
That is, $[\{\xi_{ij}\}] = [\{\xi_{ij}'\}]$, where the indices $i, j$ belong to $\{1, \dots, m\}$, 
which is the set of the indices of the open covering $\{U_i\}$.
Then there are local sections $\nu_i$ on $\mathcal L|_{U_i}$
such that 
\[
\xi_{ij} - \xi_{ij}' = \nu_i - \nu_j
\]
on $U_i\cap U_j$.
The sections $\nu_i$ determine sections on locally closed subsets $\{V_i\}$ by restriction.
%Namely, when $U_i$ does not contain a singular point of $C$, then $U_i = V_i$ and 
%$\nu_i$ obviously gives a section on $V_i$.
%When $U_i$ contains a singular point $p$ of $C$, then the restrictions of $\nu_i$
%to the branches determine sections 
%on these locally closed subsets.
By the equality $\xi_{ij} - \xi_{ij}' = \nu_i - \nu_j$ above, defining 
$\xi_i'' = \xi_i - \nu_i$, we have $\xi_{ij}'' = \xi_{ij}'$.
Then by the above argument, the classes defined by $\{\xi_{i}''\}$
and $\{\xi_{i}'\}$ give the same value of the pairing.
On the other hand, adding restrictions of sections of $\mathcal L$
 does not change the sum of residues at the singular points since $\omega_C$ is the sheaf of Rosenlicht differentials.
Thus, $\{\xi_i\}$ and $\{\xi_i'\}$ give the same value of the pairing.

Now, we prove the rest of the claims of the proposition.
Assume that $C$ has at least one singular point.
Assume $\psi$ has an $|a_1|$-th order zero or pole (when $a_1$ is non-negative or negative, 
 respectively) at some singular point $p_1$ of $C$.
Let $V_1$ be one of the disks containing $p_1$ and
 take a section $\xi_1$ of $\mathcal L|_{V_1}$
 so that it has an exactly $|a_1+1|$-th order pole or zero at $p_1$, and take all the other $\xi_i$ to be zero.
Then the paring between $\psi$ and $\{\xi_{ij}\}$ gives a nonzero value, 
 showing that the formula $\langle \psi, \{\xi_{ij}\}\rangle = \sum_{i=1}^l\sum_{j=1}^{a_i} r_{i_j} + \sum_p r_p$
 actually gives a nondegenerate pairing.
This also shows that any cohomology class in 
$H^1(C, \mathcal L)$
 can be represented by some $\{\xi_i\}$. 
When $C$ does not contain a singular point,
 we can take the open cover $\{U_i\}$ in the following way.
Namely, given any finite set (in fact, $n = 1$ suffices) of points $\{p_1, \dots, p_n\}$, we can take $\{U_i\}$ so that
 each $p_j$ is contained in the unique open subset belonging to $\{U_i\}$ (we call it $U_j$),
 and $U_j\neq U_k$ if $p_j\neq p_k$.
Then applying the same argument as above, we
 see the same conclusion holds.
  \qed

\subsection{Expression of the obstruction}\label{subsec:obstcal}
Now let us return to the calculation of the obstruction to deform $\varphi\colon C\to X$,
 where $X$ is a smooth complex surface.
We write by $\mathcal X_k$ the product space $X\times \Spec\Bbb C[t]/t^{k+1}$.
%Let $\{W_i\}$ be an open covering of $X$.
%In particular, the family $\mathcal X$ is isomorphic to the product $X\times B$ 
% over differentiable category.

Assume that we already have a $k$-th order deformation 
\[
\varphi_k\colon \mathcal C_k\to \mathcal X_k,
\]
 for some natural number $k$,
 where $\mathcal C_k$ is a flat family of curves over $\Bbb C[t]/t^{k+1}$
 whose central fiber is $C$, and 
 $\varphi_k$ restricts to $\varphi$ over $\Bbb C[t]/t$.

Let $I, J$ be sets such that $I\subset J$. 
Take an open covering $\{U_{k,i}\}_{i\in I}$ of $\mathcal C_k$ 
 and a covering $\{W_j\}_{j\in J}$ of $X$ by coordinate neighborhoods so that 
 for each $i$, the following conditions hold:
\begin{itemize}
\item The image $\varphi_k(U_{k,i})$ is contained in $W_{k, i}:=W_{i}\times\Spec\Bbb C[t]/t^{k+1}$.
\item The intersection $W_{k,i}\cap \varphi_k(U_{k,i})$ is represented by an equation
 $f_{k, i} = 0$, where $f_{k, i}$ is an analytic function on $W_{k,i}$.
\end{itemize}

On the intersection $W_{k,i}\cap W_{k,j}$ such that $U_{k, i}\cap U_{k, j}\neq\emptyset$,
 the functions $f_{k, i}$ and $f_{k, j}$ are related by
\[
f_{k, i} = g_{k, ij}f_{k, j},
\]
 where $g_{k, ij}$ is a holomorphic function on $W_{k,i}\cap W_{k,j}$
 whose reduction over $\Bbb C[t]/t$ is a non-vanishing function on $W_i\cap W_j$.
Regarding these functions as defined over $\Bbb C[t]/t^{n+2}$, we obtain the difference
\[
t^{n+1}\nu_{ij, n+1} = f_{k, i} - g_{k,ij}f_{k, j}
\]
 on $W_{k+1, i}\cap W_{k+1, j}$.
Here $\nu_{ij, n+1}$ can be regarded as a holomorphic function on $W_{i}\cap W_j$.

Assume that there is another $W_{k, l}$ containing $\varphi_{k}(U_{k, l})$,
 and define $f_{k, l}$, $g_{k, jl}$ and $g_{k, li} (= g_{k, il}^{-1})$ as above.
These determine functions $\nu_{k+1, jl}$, $\nu_{k+1, li}$ etc. on the relevant intersections.

\begin{lem}
	On $W_{i}\cap W_j\cap W_l \cap \varphi(U_i)$, 
	the identity
	\[
	\nu_{k+1, il} = \nu_{k+1, ij} + g_{0, ij}\nu_{k+1, jl}
	\]
	holds (see also \cite{KS}).
\end{lem}
\begin{rem}
	Note that 
	$W_{i}\cap W_j\cap W_l \cap \varphi(U_i)
	= W_{i}\cap W_j\cap W_l \cap \varphi(U_j)
	= W_{i}\cap W_j\cap W_l \cap \varphi(U_l)$
	holds.
	However, in general this is not equal to 
	$W_{i}\cap W_j\cap W_l \cap \varphi(C)$.
\end{rem}
\proof
By definition, we have
\[
\begin{array}{ll}
t^{n+1}\nu_{k+1, il} & = f_{k, i} - g_{k, il}f_{k, l}\\
& =  f_{k, i} - g_{k, ij}f_{k, j} + g_{k, ij}f_{k, j} -  g_{k, il}f_{k, l} \\
& = t^{k+1}\nu_{k+1, ij} + g_{k, ij}(f_{k, j} - g_{k, jl}f_{k, l}) + (g_{k, ij}g_{k, jl}-g_{k, il})f_{k, l}\\
& = t^{k+1}\nu_{k+1,ij} + t^{k+1}g_{k, ij}\nu_{k+1, jl} + (g_{k, ij}g_{k, jl}-g_{k, il})f_{k, l}.
\end{array}
\]
Also, note that the identity
\[
g_{k, il} \equiv g_{k, ij}g_{k, jl} \ \ \text{mod $t^{n+1}$}
\]
 holds.
Therefore, 
we have 
\[
(g_{k, ij}g_{k, jl}-g_{k, il})f_{k, l} = (g_{k, ij}g_{k, jl}-g_{k, il})f_{0, l} \ \  \text{mod $t^{n+2}$}.
\]
Since the identity $f_{0, l} = 0$ holds on 
$W_{i}\cap  W_j\cap W_l \cap \varphi(U_i)$, 
we have the claim.\qed\\

Note that the restriction of the 
 set of functions $\{g_{0, ij}\}$ to $\varphi(C)$ is the set of transition functions for the
normal sheaf of $\varphi$.
Thus, the lemma shows that the set of functions $\{\nu_{k+1, ij}\}$, when restricted to 
 $\varphi(C)$, behaves as a \v{C}ech 1-cocycle 
 with values in the normal sheaf of $\varphi$.
By construction, we have the following (see also \cite{KS}).
\begin{lem}\label{lem:obstcech}
The \v{C}ech 1-cocycle $\{\nu_{k+1, ij}\}$ is the obstruction cocycle 
 to deform $\varphi_k$ to a map over $\Bbb C[t]/t^{k+2}$. \qed
\end{lem}

However, it is difficult to calculate this class in this form. 
So we develop another way to express it.
Recall we fixed a function $f_{k, i}$ on $W_{k, i}$.
Regard it as a function on $W_{k+1, i}$ as above.
Then we can express it in the form
\[
f_{k, i} = f_{0, i}\exp(h_i(k+1)),
\]
 where $h_i(k+1)$ is a function on $W_{k+1, i}$
 which can have poles along $\{f_{0, i} = 0\}$ 
 (here $f_{0, i}$ is also regarded as a function on $W_{k+1, i}$), 
 and is zero when reduced over $\Bbb C[t]/t$.
Similarly, we can write $f_{k, j}$, regarded as a function on $W_{k+1, j}$,
 in the form $f_{k, j} = f_{0, j}\exp(h_j(k+1))$.
Let $h_{k+1,i}$ and $h_{k+1,j}$ be the coefficients of $t^{k+1}$ in 
$h_i(k+1)$ and $h_j(k+1)$, respectively.
These can be naturally considered as meromorphic functions on $W_i$ and $W_j$.
Then we have the following.
\begin{lem}\label{lem:residue}
The identity
\[
h_{k+1,i}-h_{k+1,j} = \frac{\nu_{k+1, ij}}{f_{0, i}} + \kappa
\]
 holds on $W_i\cap W_j$, 
 here $\kappa$ is a holomorphic function.
\end{lem}
\proof
By definition, we have
\[
t^{k+1}\nu_{k+1, ij} = f_{0, i}\exp(h_i(k+1)) - g_{k, ij}f_{0, j}\exp(h_j(k+1)).
\]
Dividing this equation by $f_{0, i}\exp(h_j(k+1))$, we have
\[
t^{k+1}\frac{\nu_{k+1, ij}}{f_{0, i}} = \exp(h_i(k+1)-h_j(k+1)) - g_{k, ij}\frac{f_{0, j}}{f_{0, i}}.
\] 
Note that since the function $\exp(h_j(k+1))$ is of the form $1+ t(\cdots)$, dividing by it does not
 affect the left hand side because the equation is defined over $\Bbb C[t]/t^{k+2}$.
Also, note that the term $g_{k, ij}\frac{f_{0, j}}{f_{0, i}}$ is holomorphic on $W_{k+1, i}\cap W_{k+1,j}$.
Therefore, the divergent terms of the function $\exp(h_i(k+1)-h_j(k+1))$ coincide with
 the divergent terms of $t^{k+1}\frac{\nu_{k+1, ij}}{f_{0, i}}$.
It follows that the function $h_i(k+1)-h_j(k+1)$ does not have a divergent term when 
 it is reduced over $\Bbb C[t]/t^{k+1}$.
Thus, the divergent terms of $h_i(k+1)-h_j(k+1)$ coincide with those of $t^{k+1}\frac{\nu_{k+1, ij}}{f_{0, i}}$.
This proves the claim.\qed\\

Let $\psi|_C$ be an element of the group $H^0(C,\varphi^*K_X)$
 and assume it is the pull back of an element $\psi$ of $H^0(X, K_X)$.
The Poincar\'e residue of the 2-form $\frac{\nu_{k+1, ij}}{f_{0, i}}\psi$ along 
$\{f_{0, i} = 0\}\cap W_i\cap W_j$
is defined by the pullback to $U_i\cap U_j$ of the 1-form $\zeta_{ij}$ on $W_i\cap W_j$
satisfying
\begin{equation}\label{eq:4}
	\frac{\nu_{k+1, ij}}{f_{0, i}}\psi = \zeta_{ij}\wedge \frac{df_{0, i}}{f_{0, i}}.
\end{equation}
From this definition, it is clear that the pullback of $\zeta_{ij}$ to $U_i\cap U_j$
coincides with the fiberwise pairing between 
$\nu_{k+1, ij}$ and $\psi|_C$ 
(recall that $\nu_{k+1, ij}$ is naturally considered as a local section 
of the normal sheaf $\mathcal N_{\varphi}$ on $ C$).
We write $\zeta_{ij}$ and its pull back to $U_i\cap U_j$ by the same letter.
These $\zeta_{ij}$ constitute a \v{C}ech 1-cocycle on $C$ with values in $\omega_C$,
 and according to Lemma \ref{lem:obstcech}, this is essentially (one of)
 the obstruction to deform $\varphi_k$.
Precisely, we have the following.
\begin{prop}\label{prop:cechobst}
If the cohomology class of $\{\zeta_{ij}\}$ in $H^1(C, \omega_C)$ 
 associated to each $\psi\in H^0(C,\varphi^*K_X)$
 vanishes, and $\varphi$ is semiregular, then the obstruction to deform $\varphi_k$ vanishes.\qed
\end{prop}

Now let us consider the meromorphic two form ${h_{k+1,i}}\psi$ on $W_i$.
We would like to define its natural analogue of the Poincar\'{e} residue along $\varphi(U_i)$
 as above. 
However, we need to take care of the facts that ${h_{k+1,i}}\psi$
 may have poles worse than logarithmic ones and also that 
 the hypersurface $\varphi(U_i)$ is singular.
In the following subsections, we study these points.

\subsection{Residues of meromorphic differential forms}
We have constructed a meromorphic two form $h_{k+1, i}\psi$ on each open subset $W_i$ of $X$.
On the intersection $W_i\cap W_j$, the difference
 $h_{k+1, i}\psi - h_{k+1, j}\psi$ equals to $\frac{\nu_{k+1, ij}}{f_{0, i}}\psi$ modulo 
 holomorphic terms.
The logarithmic two form $\frac{\nu_{k+1, ij}}{f_{0, i}}\psi$ essentially represents the obstruction
 (see Equation (\ref{eq:4})).
We will associate a \v{C}ech representative in
 the form described in Subsection \ref{subsec:paring}
 with this datum.
In the argument below, we use the same notation as in Subsection \ref{subsec:paring}.

Recall that the cohomology class of the \v{C}ech representative is determined by 
 the residues at the poles of suitable pairings.
First, we introduce a quantity which will play the role of the residues in the present context.
Let $\{U_i\}$ be the open covering of $C$ as before.
Each $U_i$ is the union of disks $U_{1, i}, \cdots, U_{k_i, i}$ attached at one point $p$,
 when $U_i$ contains a singular point of $C$.
Since $\varphi$ is an immersion, we identify $U_{a, i}$ with its image in $X$.
%Let $U_{a, i}'$ be a slightly smaller disk which is relatively compact in $U_{a, i}$.
%In particular, $U_{a, i}'$ contains the point $p$.

Fix an arbitrary Riemannian metric on $X$.
Let $\gamma$ be any smooth simple closed curve on $U_{a, i}$ which generates the fundamental group of
 the punctured disk $U_{a, i}\setminus \{p\}$.
On any open subset $\mathring U_{a, i}$ of $U_{a, i}$ away from $p$ 
 (that is, the closure of $\mathring U_{a, i}$ does not contain $p$), the exponential map 
 gives a diffeomorphism from a small neighborhood of the zero section of the normal bundle
 of $\mathring U_{a, i}$ to a tubular neighborhood of $\mathring U_{a, i}$ in $X$.

Let $S_{\delta}\mathring U_{a, i}$ be the circle bundle of radius $\delta$ in the normal bundle.
Here, $\delta$ is a small positive number.
By taking $\delta$ small enough, we can assume the image of $S_{\delta}\mathring U_{a, i}$ by the
 exponential map is disjoint from $\varphi(U_i)$.
Then, the meromorphic 2-form $h_{k+1, i}\psi$ is pulled back to a smooth 2-form on 
 $S_{\delta}\mathring U_{a, i}$.
We write the pull back by the same notation $h_{k+1, i}\psi$ to save letters.
Let $T_{\gamma}$ 
 be the union of the fibers of $S_{\delta}\mathring U_{a, i}$ over $\gamma$, which is diffeormorphic to 
 the two dimensional torus.
Note that the fiber of $T_{\gamma}$ has a natural (positive) orientation 
 induced from the natural complex orientation of the normal bundle of $\mathring U_{a, i}$.
This and the orientation of $\gamma$ determine an orientation of $T_{\gamma}$ by 
 ordering a basis of a tangent space as \{base, fiber\}.
This is chosen so that it will be compatible with Equation (\ref{eq:4}) in the following calculation.

\begin{lem}\label{lem:resint}
The integration of $h_{k+1, i}\psi$ over $T_{\gamma}$ does not depend on the metric, $\delta$, or $\gamma$.
Therefore, it only depends on $p$ and $a$.
\end{lem}
\proof
This follows from the Stokes' theorem since $h_{k+1, i}\psi$ is a closed 2-form.\qed

\begin{defn}\label{def:residualvalue}
When $\gamma$ is positively oriented with respect to 
 the complex orientation of $U_{a, i}$, 
 then we write the value obtained in Lemma \ref{lem:resint} by $r(p, a)$.
\end{defn}

Let $C_l$ be any irreducible component of $C$.
Let $\pi_l\colon \widetilde C_l\to C_l$ be its normalization.
Let $s\widetilde C_l$ be the set of points on $\widetilde C_l$ which are mapped to 
 singular points of $C$ by $\pi_l$.
Also, let $\mathring C_l$ be an open subset of $\widetilde C_l$ obtained by deleting 
 small disks around $s\widetilde C_l$.
Note that $\mathring C_l$ can also be seen as an open subset of $C_l$, and we assume $\mathring C_l$ is covered by
 open subsets of the form $U_{a, i}$. 
On such an open subset $\mathring C_l$, the above construction glues and gives a circle bundle
 $S_{\delta}\mathring C_l$
 of radius $\delta$ inside the normal bundle of $\mathring C_l$.

On each intersection $U_{a,i}\cap \mathring C_l$, the 2-form $h_{k+1, i}\psi$ is pulled back to 
 $S_{\delta}\mathring C_l|_{U_{a, i}\cap \mathring C_l}$
 and gives a smooth 2-form, which again we write by the same notation.
Then since the fiber integration commutes with the exterior derivative,
 that of $h_{k+1, i}\psi$ gives a closed 1-form on $U_{a, i}\cap \mathring C_l$,
 which we write by $\int_{\delta} h_{k+1, i}\psi$.

When $U_{a, i}$ does not contain a singular point of $C$ (so that $U_i = U_{a, i}$ is really a disk),
 then the form $\int_{\delta} h_{k+1, i}\psi$
 is defined on the whole $U_{a, i}$.
When $U_{a, i}$ contains a singular point $p$,
 $\int_{\delta} h_{k+1, i}\psi$ is defined only on $U_{a, i}\cap \mathring C_l$.
But taking $\mathring C_l$ larger by taking $\delta$ smaller, 
 we can assume that the intersection $U_i\cap U_j\cap C_l$ is contained in $\mathring C_l$
 for any distinct $i$ and $j$.
Since $\int_{\delta} h_{k+1, i}\psi$ is a closed form,
 the value of the integration of it along a contour which circles $p$ once in the positive direction
 does not depend on the contour and given by $r(a, p)$ by Lemma \ref{lem:resint}.
% we can take a tubular neighborhood $U$ of $U_{a, i}\setminus\{p\}$
% and a map $\pi\colon U\to U_{a, i}\setminus\{p\}$
% such that $\pi$ is a holomorphic disk bundle.
% let $\rho$ be a smooth cut off function on $U_{a, i}$ such that
%\begin{itemize}
%\item $0\leq \rho \leq \delta$,
%\item $\rho(p) = 0$ and $\rho \equiv \delta$ on $U_{a, i}\cap \mathring C_l$, and
%\item $\rho(x) \neq 0$ for $x\neq p$.
%\end{itemize}
%By shrinking each $U_i$ slightly if necessary, we can assume that the support of $d\rho$
% is disjoint from $U_j$, for all $j\neq i$.
%We will assume this throughout.
%Taking such a function $\rho$ suitably and considering a circle bundle over $C_l\setminus \{p\}\cap U_{a, i}$
% whose fiber over $x$ is the circle of radius $\rho(x)$, the circle bundle embeds in 
% $X\setminus \varphi(U_i)$ by the exponential map.

%Thus, the 2-form $h_{k+1, i}\psi$ is pulled back to a smooth 2-form
% on the circle bundle (recall that $h_{k+1, i}\psi$ may be singular only on $\varphi(U_i)$).
%The fiber integration of this form gives a closed 1-form on $U_{a, i}\setminus\{p\}$
% which extends $\int_{\delta} h_{k+1, i}\psi$, which we still write by 
% $\int_{\delta} h_{k+1, i}\psi$ for notational simplicity.

Now consider the intersection $U_i\cap U_j\cap C_l$.
According to Lemma \ref{lem:residue}, we have
\[
\int_{\delta} h_{k+1, i}\psi - \int_{\delta} h_{k+1, j}\psi
 = \int_{\delta} \left(\frac{\nu_{k+1, ij}}{f_{0, i}} + \kappa\right)\psi.
\]
These naturally give a $\mathcal C^1(\widetilde C_l)$-valued \v{C}ech 1-cocycle on 
 $\widetilde C_l$ associated with a covering of $\widetilde C_l$ induced by $\{U_i\}$.
Here $\mathcal C^1(\widetilde C_l)$ is the sheaf of complex valued $C^{\infty}$ closed 1-forms
 on $\widetilde C_l$.
We will compare the cohomology class defined by this cocycle with the ones described in 
 Subsection \ref{subsec:paring}.

\subsubsection{Cohomology classes defined by closed 1-forms and the pairing}

In general,
 let $S$ be a compact nonsingular Riemann surface with finite number of 
 marked points $\{p_i\}_{i=1}^b$, $b\geq 1$.
Take an open covering $\{U_j\}$ of $S$ so that each $p_i$ is contained in 
 a unique open subset belonging to  $\{U_j\}$, which we write by $U_i$.
We assume $U_i\neq U_j$ if $p_i\neq p_j$.
When $U_j$ does not contain a marked point, then we associate with it a closed 1-form $\mu_j$ on it.
When $U_j$ contains a marked point $p_j$, then we associate with it a closed 1-form $\mu_j$ on
 $U_j\setminus B(p_j)$, where $B(p_j)$ is a small disk around $p_j$.
We take each disk $B(p_j)$ small enough so that the intersection $U_i\cap U_j$ for distinct $i$ and $j$
 does not intersect any such disk. 
Note that since $\mu_i$ is closed, the contour integral along a closed path in $U_j$ is zero when 
 $U_j$ does not contain a marked point, and the integral along a closed path in $U_j\setminus B(p_j)$
 depends only on the homotopy class of the path when $U_j$ contains the marked point $p_j$.
When a path circles around $p_j$ once in the positive direction, we write by $r(p_j)$ the 
 value of the contour integral of $\mu_j$ along the path.

Such a family of closed 1-forms $\{\mu_i\}$ determines a \v{C}ech 1-cocycle 
 with values in $\mathcal C^1(S)$ by associating $\mu_{ij}:=\mu_i-\mu_j$
 with the intersection $U_i\cap U_j$.
There is a standard resolution of the sheaf $\mathcal C^1( S)$, 
\[
0\to \mathcal C^1(S)\to \mathcal A^1(S)\to \mathcal A^2(S)\to 0
\]
 by flabby sheaves.
Here $\mathcal A^i(S)$ is the sheaf of complex valued $C^{\infty}$ $i$-forms
 on $S$.
In particular, the cohomology group $H^1(S, \mathcal C^1(S))$
 is isomorphic to $H^2(S, \Bbb C) \cong H^1(S, \omega_{S})$, 
 where $\omega_{S}$ is the canonical sheaf of $S$.

Since $H^2(S, \Bbb C)$ is dual to $H^0(S, \Bbb C)$, the \v{C}ech 1-cocycle 
 $\{\mu_{ij}\}$ should give an element of the dual of $H^0(S, \Bbb C)$ in a natural manner.
This is an analogue of Proposition \ref{prop:pair}.
We need to associate a scalar with a constant function $\alpha\in H^0(S, \Bbb C)$ in a natural way.
Namely, we associate the value
\[
\langle \{\mu_{ij}\}, \alpha\rangle = \alpha\sum_i  r(p_i)
\] 
 with it.
%First, consider an open subset $U_j$ and a function $s_j$ with a compact support on it,
% such that it is constant in a small neighborhood of a singular point (if $U_j$ contains it).
%If $U_j$ does not contain a singular point, then 
% we associate the value 
%\[
%\int_{U_j} ds_j\wedge \mu_j 
%\]
% with it, which is zero by the Stokes' theorem.
%
%Suppose that $U_i$ contains a singular point $p_i$.
%Let $B_{\varepsilon_n}(p_i)$, $n\in \Bbb N$ be a sequence of disks around $p_i$
% whose radii $\varepsilon_n$ converges to zero as $n$ goes to infinity.
%Then, to a test function $s_i$ as above, we associate 
% the value
%\[
%\lim_{n\to \infty}\int_{U_i\setminus B_{\varepsilon_n}(p_i)} ds_i\wedge \mu_i.
%\]
%By standard calculation, it is not difficult to see that this converges to 
% the value $s_i(p_i)r(p_i)$.
% 
%Now consider a constant function $\alpha\in H^0(S, \Bbb C)$.
%Take a partition of unity $\rho_j$ associated with the covering $\{U_j\}$
% (after refining $\{U_j\}$ if necessary).
%Then we associate the value
%\[
%\langle \{\mu_{ij}\}, \alpha\rangle = \alpha\sum_i \rho_i(p_i)r(p_i)
%\] 
% to it.

\begin{prop}\label{prop:pair2}
This pairing between the classes $\{\mu_{ij}\}\in H^2(S, \Bbb C)$ and $\alpha\in H^0(S, \Bbb C)$
 is well-defined and non-trivial.
\end{prop}
\proof
We  need to check that the value $\langle \{\mu_{ij}\}, \alpha\rangle$
 does not depend on the choices of the 1-forms $\{\mu_i\}$
 with $\{\mu_{ij}\}$ fixed, nor the representative $\{\mu_{ij}\}$ of the cohomology class.

%The independence from the choice of the partition of unity is straightforward and we omit it.
Given $\{\mu_{ij}\}$, assume we have another set of singular 1-forms $\{\mu'_i\}$
 satisfying $\mu'_i-\mu_j' = \mu_{ij}$.
As in the proof of Proposition \ref{prop:pair}, the set of differences $\{\mu_i-\mu'_i\}$
 gives a  global closed 1-form $\mu$ on $\mathring S:= S\setminus (\cup_j B(p_j))$.
If $\{\gamma_i\}$ is the set of disjoint contours around $p_i$ encircling $p_i$ once in the 
 positive direction, then by the Stokes' theorem, we have
\[
\sum_i\int_{\gamma_i}\mu = 0.
\]
Therefore, the pairing between $\mu$ and $\alpha$ is zero.
Thus, the pairing does not depend on the choice of  $\{\mu_i\}$ with $\{\mu_{ij}\}$ fixed. 
 
Suppose there is another representative $\{\mu'_{ij}\}$ of the cohomology class $[\{\mu_{ij}\}]$ 
 in $H^1(S, \mathcal C^1(S))$.
Then there is a set of smooth closed 1-forms $\{\nu_i\}$ on $U_i$ such that 
 $\mu_{ij} - \mu'_{ij} = \nu_i-\nu_j$ on $U_i\cap U_j$.
Then defining $\mu''_i: = \mu_i-\nu_i$, we have $\mu''_{ij} = \mu'_{ij}$.
From this, it is clear that the pairing between $\{\mu_{ij}\}$ and $\alpha$ is equal to 
 that between $\{\mu'_{ij}\}$ and $\alpha$.
 
The proof of the fact that the paring is non-trivial is given by
 an obvious modification of the proof of Proposition \ref{prop:pair}.\qed\\

By definition of the pairing above, we see the following.

\begin{cor}\label{cor:resbasis}
The class defined by $\{\mu_{ij}\}$ can be identified with 
 $\sum_i r(p_i)$ using the basis of $H^2(S, \Bbb C)$
 dual to the standard generator of $H^0(S, \Bbb Z)\subset H^0(S, \Bbb C)$.\qed
\end{cor}

\subsubsection{The main theorem}
Applying the above construction 
 to our situation,
 we have the following.
\begin{prop}
	The cohomology class of $H^1(\widetilde C_l, \mathcal C^1(\widetilde C_l)) 
	\cong H^2(\widetilde C_l, \Bbb C)$ defined by 
	the \v{C}ech 1-cocycle
	$\left\{ \int_{\delta} h_{k+1, i}\psi - \int_{\delta} h_{k+1, j}\psi    \right\}$
	does not depend on $\delta$ and identified with
	$\sum_{a,p} r(a, p)$ in the sense of Corollary \ref{cor:resbasis}.
	Here $p$ runs through the singular points of $C$ on $C_l$ and $a$ indexes the 
	local branches of $C$ at $p$ contained in $C_l$.  \qed
\end{prop}

In particular, if the limit $\lim_{\delta\to 0}\left\{ \int_{\delta} h_{k+1, i}\psi - \int_{\delta} h_{k+1, j}\psi    \right\}$
 exists, then it will give the same cohomology class as 
 $\left\{ \int_{\delta} h_{k+1, i}\psi - \int_{\delta} h_{k+1, j}\psi    \right\}$.
This is in fact the case by 
 the equality $\int_{\delta} h_{k+1, i}\psi - \int_{\delta} h_{k+1, j}\psi
 = \int_{\delta} \left(\frac{\nu_{k+1, ij}}{f_{0, i}} + \kappa\right)\psi$,
 since the integral $\int_{\delta} \left(\frac{\nu_{k+1, ij}}{f_{0, i}} + \kappa\right)\psi$
 converges to the usual Poincar\'e residue $\zeta_{ij}$ (see Equation  (\ref{eq:4}))
 as $\delta$ goes to zero.
 
The set of local 1-forms $\left\{ \int_{\delta} h_{k+1, i}\psi - \int_{\delta} h_{k+1, j}\psi    \right\}$
determines a class of $H^1(\widetilde C_l, \mathcal C^1(\widetilde C_l))\cong H^1(\widetilde C_l, \omega_{\widetilde C_l})$
 for each component $\widetilde C_l$,
 while the set of holomorphic 1-forms $\{\zeta_{ij}\}$ determines a class of $H^1(C, \omega_C)$.
Although these groups are a priori different, the restriction of a class in $H^1(C, \omega_C)$
to the component $C_l$ is naturally identified with the dual of 
$H^0(C_l, \mathcal O_{C_l})\cong H^0(C_l, \Bbb C)\cong H^0(\widetilde C_l, \Bbb C)$, 
using the description in Subsection \ref{subsec:paring}.
By Proposition \ref{prop:pair}, the class $\{\zeta_{ij}\}$ can be represented by 
local meromorphic 1-forms on open subsets $U_{i, a}$.
Regarding these 1-forms as closed singular $C^{\infty}$ 1-forms, it is clear that this identification 
is compatible with the one $H^1(\widetilde C_l, \mathcal C^1(\widetilde C_l))\cong H^0(\widetilde C_l, \Bbb C)^{\vee}$
described above.
Thus, the set of local 1-forms $\left\{ \int_{\delta} h_{k+1, i}\psi - \int_{\delta} h_{k+1, j}\psi    \right\}$
also determines a class of $H^1(C, \omega_C)$.
Therefore, we have the following.
\begin{cor}
The cohomology classes of $H^1(C, \omega_C)$
 determined by $\left\{ \int_{\delta} h_{k+1, i}\psi - \int_{\delta} h_{k+1, j}\psi    \right\}$
 and $\{\zeta_{ij}\}$ coincide. \qed
\end{cor}

By this corollary, to prove the vanishing of the obstruction class $\{\zeta_{ij}\}$, 
 it suffices to see the vanishing of the class 
 $\left\{ \int_{\delta} h_{k+1, i}\psi - \int_{\delta} h_{k+1, j}\psi    \right\}$
 for all $\psi\in H^0(X, K_X)$, when semiregularity holds.
Recall that this is the same as the vanishing of the scalar $\sum_{a,p} r(a, p)$.
Thus, if the sum 
 $\sum_a r(a, p)$ vanishes  for each fixed $p$, the obstruction vanishes, too.
Here, the sum runs over all the branches of $C$ at $p$.
Recall that the scalar $r(a, p)$ is the value of the integration of the two form $h_{k+1, i}\psi$
 along the restriction of the circle bundle $S_{\delta} \mathring C_l$ 
 over a loop $\gamma_{a, p}$ encircling the singular point $p$
 (here $p$ is contained in the open subset $U_{a, i}\subset C_l$).
The circle bundle is embedded in $X$ by the exponential map.

Now consider the boundary $\partial B\cong S^3$ of a small ball in $X$ around $p$.
Let $D_{\delta}\gamma_{a, p}$ be the disk bundle over $\gamma_{a, p}$ whose boundary
 is $S_{\delta} \mathring C_l|_{\gamma_{a, p}}$.
Topologically it is a solid torus.
By suitably isotoping, these disk bundles can be disjointly embedded in $\partial B$ 
 for all branches indexed by $a$.
Now by the Stokes' theorem, we have
\[
\sum_a r(a, p) = \sum_a \int_{S_{\delta} \mathring C_l|_{\gamma_{a, p}}} h_{k+1, i}\psi
 = -\int_{\partial B\setminus \cup_a D_{\delta}\gamma_{a, p}} d(h_{k+1, i}\psi)
 = 0.
\]

Thus, we have proved the following theorem.
Let $\varphi\colon C\to X$ be a map from a reduced complete curve to a smooth complex surface
 which is an immersion.

\begin{thm}\label{thm:main}
If the map $\varphi\colon C\to X$ is 
 semiregular, then it is unobstructed in the sense that 
 any first order deformation
 can be extended to any higher order.\qed
\end{thm}

This enables us to construct large number (often infinitely many) of curves
with specific geometric properties from a small number of data.

\begin{rem}\label{rem:locality}
Note that the calculation of the obstruction to deform $\varphi$ is local.
Namely, whether $\varphi$ is obstructed or not is determined by 
 the data of a neighborhood of $\varphi(C)$ in $X$.
In particular, we may replace $X$ by a suitable open subset $X'$ of $X$.
Even if $\varphi$ is not semiregular as a map $\varphi\colon C\to X$
 (for example, when $H^{0}(X, K_X)$ vanishes but $H^{0}(C, \varphi^*K_X)$
 does not), when some $\varphi\colon C\to X'$, with the target suitably replaced,
 is semiregular, then the original $\varphi$ is also unobstructed.

%By the same reason, the curve $C$ can be more general than the argument above, where
% we assumed $C$ to be compact (and reduced).
%Namely, $C$ can contain some punctures and boundary components, and if we combine with the method in \cite{N1},
% we can use Theorem \ref{thm:main} to produce many disks with Lagrangian boundary conditions
% in various situations.
\end{rem}

\subsection{Higher dimensional case}\label{subsec:higherdim}
We can apply the above argument to the case of 
	a hypersurface $\varphi\colon Z\to X$, when $Z$ is smooth
	and $\varphi$ is an immersion.
In this case, the semiregularity map is the natural map
\[
H^{n-2}(X, K_X)\to H^{n-2}(Z, \varphi^*K_X).
\]
Here $n$ is the dimension of $X$.
\begin{thm}
If this is surjective, then the map $\varphi$ is unobstructed.
\end{thm}
\proof
The proof in the previous subsection applies with obvious modifications, 
 so we only give a sketch.
In this case, the covering $\{U_i\}$ of $Z$ consists of smooth open subsets, 
 and the integration along the fibers of the circle bundle introduced after 
 Definition \ref{def:residualvalue}
 produces a closed $(2n-3)$-form on $U_i$.
Therefore, the closed $(2n-3)$-form valued \v{C}ech 1-cocycle defined by the 
 differences of such forms on $\{U_i\}$ is cohomologically trivial by definition.
This cohomology class, when integrated over the fundamental class of $Z$,
 converges to the obstruction class 
 coupled with an element of the dual space $H^{n-2}(Z, \varphi^*K_X)$
 when the radius of the circle bundle goes to zero.
Therefore, the obstruction vanishes.\qed

\begin{rem}
This result was proved in \cite{I} when $Z$ and $X$ are K\"ahler.
There, the vanishing of non-curvilinear obstructions was also proved. 
\end{rem}

\section{Deformation of maps into Calabi-Yau surfaces}\label{sec:4}
Let $\varphi\colon C\to X$ be a map from a reduced connected curve to a smooth surface
 which is an immersion as before.
Theorem \ref{thm:main} applies to various situations.
The most simple case is 
 when the canonical sheaf of $X$ is trivial, because in this case \emph{any} $\varphi$ is
 semiregular, and thus unobstructed.
Before applying this to actual problems,
 we introduce the following terminology.

Let us write $C = \cup C_i$, where $C_i$ is an irreducible component of $C$.
Let $\widetilde C = \cup \widetilde C_i$ be the normalization of $C$.
Here $\widetilde C_i$ is the normalization of $C_i$.

\begin{defn}
Let $C'$ be a nodal curve.
We call it an \emph{unchaining} of $C$ if the
 natural map $p\colon \widetilde C\to C$ factors through $C'$.
\end{defn}

Let $C'$ be an unchaining of $C$ and assume $C$ is nodal.
Then the map $\varphi\colon C\to X$ gives a natural map 
 $\varphi'\colon C'\to X$ which is also an immersion.
In general it is more difficult to deform $\varphi'$ than $\varphi$.
However, as noted above, when the canonical sheaf $K_X$ 
 is trivial, they are both unobstructed.

This simple claim can be used to produce a lot of curves of prescribed genera on Calabi-Yau surfaces.
For example, let us consider the following result, see \cite{H}, Chapter 13, Proposition 2.4.
\begin{prop}\label{prop:K3nodalcurve}
Let $X$ be a K3 surface over an algebraically closed field of characteristic zero.
Suppose there exists an integral nodal rational curve $C$ in $X$ of arithmetic genus $g>0$.
Then there is a one-dimensional family of nodal elliptic curves in the linear system $|\mathcal O(C)|$. \qed
\end{prop}
In \cite{H}, this result is proved by comparing the dimensions of $|\mathcal O(C)|$
 and the moduli space of stable curves.
Using Theorem \ref{thm:main}, we can prove a stronger result with a simpler proof. 
See also \cite[Theorem 1.1]{D} for the case where $C$ is ample and irreducible.
\begin{prop}\label{prop:embrational}
Let  
 $X$ be a complex K3 surface
 and $C$ be a connected, embedded, reduced nodal rational curve of arithmetic genus $g$ in $X$
 ($C$ need not be irreducible).
Then there is a $k$ dimensional family of nodal curves of geometric genus $k$
 for each $0\leq k\leq g$.
\end{prop}
\proof 
Consider the curve $C$
 as an embedding $\varphi\colon C\to X$.
Take a connected unchaining $C'$ of $C$ so that the arithmetic genus of 
 $C'$ is one.
By Theorem \ref{thm:main}, the induced map $\varphi'\colon C'\to X$
 is unobstructed.
It is easy to see that the dimension of $H^0(C', \nu_{\varphi'})$
 is one.
Thus, the map $\varphi'$ has a one dimensional nontrivial deformation.
Since a complex K3 surface cannot be dominated by a family of rational curves, 
 the image of a deformed map of $\varphi'$ must be an elliptic curve.
 
Next, consider a connected unchaining $C''$ of $C$ so that the arithmetic genus of 
 $C''$ is two.
As before, the induced map $\varphi''\colon C''\to X$
 is unobstructed, and one sees that the dimension of $H^0(C'', \nu_{\varphi''})$
 is two.
Since there are only one dimensional families of elliptic curves deforming $C$ by the above paragraph, 
 the image of a general deformation of $\varphi''$ has geometric genus two.
Repeating the argument in the same way, one obtains the proposition.\qed\\

Pursuing the idea of the proof of this proposition, we can construct many curves
 with controlled geometry from a few number of curves.
After this paper was posted on arXiv, \cite{CGL} proved strong existence theorem of 
 infinitely many integral curves of given geometric genus in any projective K3 surface.
Although we can only produce such curves in a Zariski open subset in the moduli space
 of K3 surfaces, we give the construction since our method is different and constructive, and
 it gives an example of controlled deformations of non-reduced curves, which cannot
 be done with the previously known techniques.

Let $\alpha\in H^2(X, \Bbb Z)$ be the class of $X$ which gives the polarization.
In the paper \cite{C}, Chen proved that for every positive integer $k$, 
 there is a Zariski open subset $U_k$ of $\mathcal M$ whose member contains
 (finite number of) irreducible nodal rational curves of the class $k\alpha$.
Moreover, for different $k$ and $k'$, the intersection 
 $U_k\cap U_{k'}$ contains a Zariski open subset 
 on which the intersection of the rational curves of degrees $k$ and $k'$
 are ordinary nodes.

\begin{figure}[h]
	\includegraphics[height=4cm]{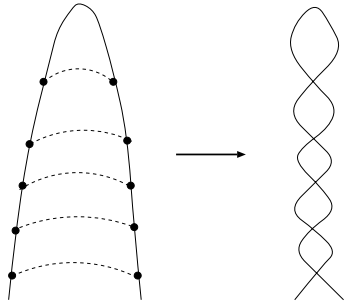}
	\caption{{A nodal rational curve $C$ in a K3 surface $X$ (the figure on the right)
	               is seen as the image of a map from its normalization $\tilde C$.
	               The dotted lines indicate the conductors.
                   In this figure, the curve looks irreducible, but in general the curve need not be irreducible.}}\label{fig:1}
\end{figure}

Using this, it is clear that there is a Zariski open subset of $\mathcal M$ whose
 member contains an embedded (not necessarily irreducible) nodal rational curve $C$ of arithmetic genus
 at least two.
%Let $\varphi\colon C\to X$ be the inclusion.
Consider the normalization $\pi\colon \tilde C\to C$ of it, see Figure \ref{fig:1}.

%As in the proof of Proposition \ref{prop:embrational}, 
% this map $\varphi_1$ deforms and gives a one parameter family of nodal curves of 
% geometric genus one.

\begin{figure}[h]
	\includegraphics[height=8cm]{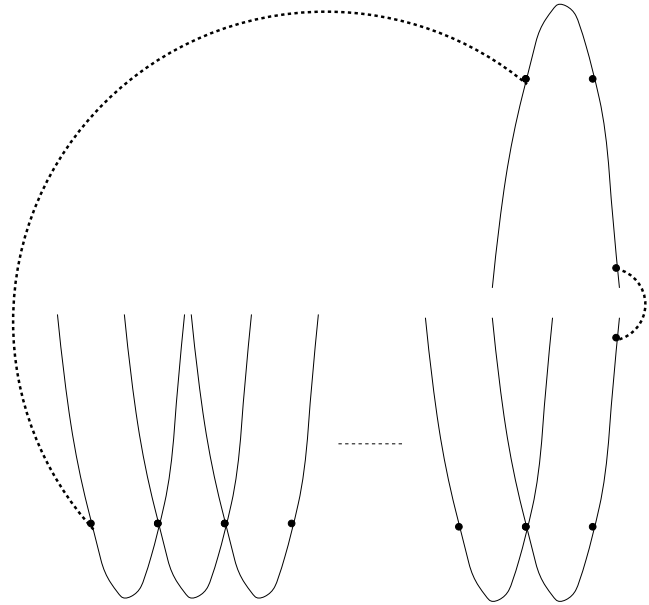}
	\caption{{A chain of rational curves of arithmetic genus one.
	Each component is a copy of $\tilde C$ in Figure \ref{fig:1}.
    The dotted lines indicate that the connected points are glued.}}\label{fig:1.5}
\end{figure}

Next, using the normalized curve, we construct a chain of rational curves of arithmetic genus one, 
 see Figure \ref{fig:1.5}. 

%\begin{figure}[h]
%	\includegraphics[height=4cm]{rationalfig2.eps}
%	\caption{The points $P$ and $Q$ are identified.
%	               The result is a nodal rational curve of arithmetic genus one.}\label{fig:2}
%\end{figure}

We write this curve as $C_k$ when its number of components is $k\geq 2$.
The map $\pi$ naturally gives a map $\varphi_k\colon C_k\to X$.
As before, this curve deforms and the image is a one parameter family of immersed curves of
 geometric genus one.
Note that under the deformation, all the nodes of the domain curve must be smoothed
 in order to acquire non-zero genus.
In particular, the domain curve of the deformed map and its image as well are irreducible.
Also, the map $\varphi_k$ cannot be a covering of a map of lower degree, since the domain is not
 a covering of another curve. 
This proves the theorem in the case $g = 1$.\\  

\begin{figure}[h]
	\includegraphics[height=4cm]{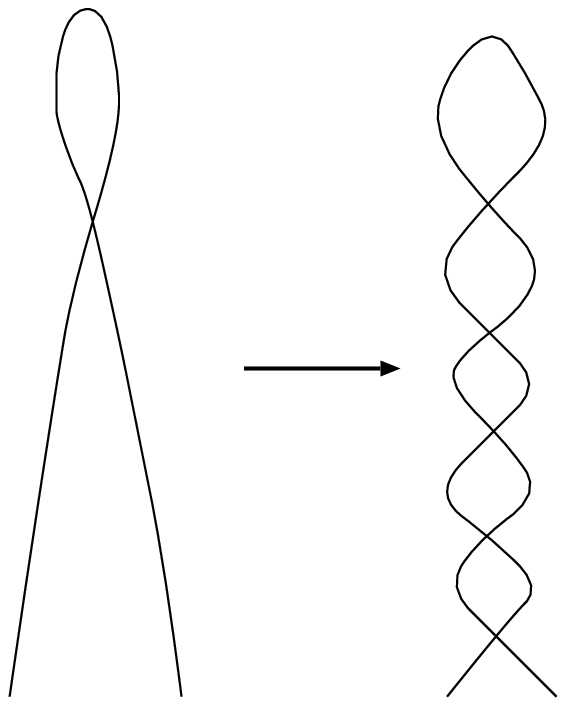}
	\caption{}\label{fig:3}
\end{figure}

Next, consider the case $g = 2$.
Recall that we took a nodal rational curve of arithmetic genus at least two.
Therefore, there is an unchaining $C'$ of it which is of arithmetic genus one
 such that the map $\pi$ above splits through it, 
 see Figure \ref{fig:3}.

\begin{figure}[h]
	\includegraphics[height=4cm]{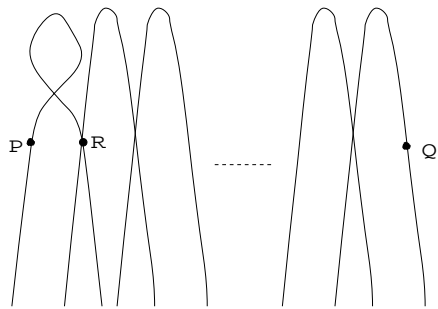}
	\caption{{The points $P$ and $R$ are the inverse images of the node of $C$ under the partial normalization.
                   The points $P$ and $Q$ are identified so that 
                  the resulting curve has arithmetic genus two.}}\label{fig:4}
\end{figure}

Then construct a chain of rational curves similar to the case $g = 1$ above using 
 a piece of $C'$ and $k-1$ pieces of $\tilde C$, see Figure \ref{fig:4}.
We obtain a prestable rational curve of arithmetic genus two. 

We write this curve as $C_k'$.
As before, the map $\pi$ induces a map $\varphi_k'\colon C_k'\to X$, 
 and it has two dimensional deformations.
As in the proof of Proposition \ref{prop:embrational}, 
 a general member of the deformed curves must be
 of geometric genus two.
This proves the theorem for $g = 2$.
For curves with higher genera, we can use the similar construction, replacing 
 more sheets in Figure \ref{fig:4} by $C'$. \qed

\section{Deformation of maps into general surfaces}
\subsection{Geometric criterion for semiregularity}
Now let us consider the application of Theorem \ref{thm:main} to more general situations.
Let $\varphi\colon C\to X$ be a map from a reduced connected nodal curve to a smooth surface 
which is an immersion.
We assume $\varphi(C)$ is a reduced nodal curve.
%The curve $C$ need not be connected.
%Let $C'$ be an unchaining of $C$ and $\varphi'\colon C'\to X$ be a natural map
% induced from $\varphi$.
% $\varphi'$ is also an immersion.
Let $\iota\colon \varphi(C)\to X$ be the inclusion.
Let $p\colon C\to \varphi(C)$ be the natural map.

Consider the exact sequence
\[
0\to \iota^*\mathcal K_X\to p_*\varphi^*\mathcal K_X\to \mathcal Q\to 0
\]
of sheaves on $\varphi(C)$.
Here $\mathcal K_X$ is the canonical sheaf of $X$ and 
$\mathcal Q$ is defined by this sequence.
Let $P = \{p_i\}$ be the set of nodes of $\varphi(C)$ whose inverse image by $p$ consists of two points.
Then the sheaf $\mathcal Q$ is isomorphic to the direct sum of skyscraper sheaves 
$\oplus_i\Bbb C_{p_i}$.
The associated cohomology exact sequence is
\[
0\to H^0(\varphi(C), \iota^*\mathcal K_X)\to H^0(\varphi(C), p_*\varphi^*\mathcal K_X)
\to \oplus_i\Bbb C_{p_i}\to H^1(\varphi(C), \iota^*\mathcal K_X)\to H^1(\varphi(C), p_*\varphi^*\mathcal K_X)\to 0.
\] 

By the Leray spectral sequence, we have
$H^i(\varphi(C), p_*\varphi^*\mathcal K_X) \cong H^i(C, \varphi^*\mathcal K_X)$, for $i = 0, 1$.
It follows that if the curve $\varphi(C)$ is semiregular in the classical sense, then the map $\varphi$
is semiregular if and only if the map 
$a\colon H^0(\varphi(C), \iota^*\mathcal K_X)\to H^0(\varphi(C), p_*\varphi^*\mathcal K_X)$ is surjective.
By the exact sequence, this condition is 
equivalent to the claim that the map
$b\colon \oplus_i\Bbb C_{p_i}\to H^1(\varphi(C), \iota^*\mathcal K_X)$
is injective.
Now we can prove the following.
\begin{thm}\label{thm:CB}
	Assume that $\varphi(C)$ is semiregular in the classical sense.
	Then the map $\varphi$ is semiregular if and only if for each $p_i\in P$, there is a first order deformation of 
	$\varphi(C)$ which smoothes $p_i$, but does not smooth the other nodes of $P$. 
\end{thm}
\proof
Note that the dual of $H^1(\varphi(C), \iota^*\mathcal K_X)$
is isomorphic to $H^0(\varphi(C), \mathcal N_{\iota})$
and the dual of $H^1(C, \varphi^*\mathcal K_X)$
is isomorphic to $H^0(C, \mathcal N_{\varphi})$.
Here $\mathcal N_{\iota}$ is the normal sheaf of the map $\iota$.
Taking the dual of the above exact sequence, the map $b$ is injective if and only if
the map 
\[
b^{\vee}\colon H^0(\varphi(C), \mathcal N_{\iota})\to \oplus_i\Bbb C_{p_i}
\]
is surjective.
By definition, elements of $H^0(C, \mathcal N_{\varphi})$ corresponds to first order deformations
of $\varphi(C)$ which do not smooth the nodes in $P$.
Then the image of $b^{\vee}$
contains the summand $\Bbb C_{p_i}$ if and only if there is a first order deformation of $\varphi(C)$
(that is, an element of $H^0(\varphi(C), \mathcal N_{\iota})$)
which smoothes the node $p_i$ but does not smooth the other nodes of $P$.
This proves the claim. \qed\\

Combining with Theorem \ref{thm:main}, we have the following.
\begin{cor}\label{cor:CB}
	Assume that $\varphi(C)$ is semiregular in the classical sense.
	Then the map $\varphi$ is unobstructed if the condition of Theorem \ref{thm:CB} is satisfied.\qed
\end{cor}

For applications, it will be convenient to write this in the following geometric form.
Consider the exact sequence
\[
0\to \mathcal O_X\to \mathcal O_X(\varphi(C))\to \mathcal N_{\iota}\to 0
\]
 of sheaves on $X$ and the associated cohomology sequence
\[
0 \to H^0(X, \mathcal O_X)\to H^0(X, \mathcal O_X(\varphi(C)))\to H^0(\varphi(C), \mathcal N_{\iota})
 \to H^1(X, \mathcal O_X)\to \cdots.
\]
Let $V$ be the image of the map $H^0(\varphi(C), \mathcal N_{\iota})
 \to H^1(X, \mathcal O_X)$.
Since we are working in the analytic category, we have the exponential exact sequence
\[
0\to \Bbb Z\to \mathcal O_X\to \mathcal O_X^{\ast}\to 0
\]
 of sheaves on $X$.
Let $\bar V$ be the image of $V$ in $Pic^0(X) = H^1(X, \mathcal O_X^{\ast})$.

\begin{cor}\label{cor:geomCB}
	In the situation of Theorem \ref{thm:CB}, the map $\varphi$ is unobstructed
	if for each $p_i\in P$, there is an effective divisor $D$ such that $\mathcal O_X(\varphi(C)-D)\in \bar V$
	and avoids $p_i$ but passes through all points in $P\setminus\{p_i\}$. 
\end{cor}
\proof
Fix a coordinate system $\{U_i\}$ on $X$ and let $\{f_i\}$ be the local defining functions of $\varphi(C)$
 on each $U_i$.
On those $U_i$ which do not intersect $\varphi(C)$, we take $f_i = 1$.
Then the meromorphic function $\frac{1}{f_i}$ on $U_i$ gives a local frame for the invertible sheaf 
 $\mathcal O_X(\varphi(C))$
 on $U_i$.
The restriction of it to $U_i\cap \varphi(C)$ gives a local frame for the normal sheaf $\mathcal N_{\iota}$.
Using this notation, a section of $\mathcal N_{\iota}$ is given by a set of 
 holomorphic functions $\{g_i\}$, where $g_i$ is defined on $U_i\cap \varphi(C)$,
 such that $\frac{g_i}{g_j} = \frac{f_i}{f_j}$ for any $i\neq j$.
Then the map $H^0(\varphi(C), \mathcal N_{\iota})
 \to H^1(X, \mathcal O_X)$ is given as follows.
Namely, extend $g_i$ to a holomorphic function $\tilde g_i$ on $U_i$.
Again, on those $U_i$ which do not intersect $\varphi(C)$, we take $\tilde g_i = 1$.
The difference $\frac{\tilde g_i}{f_i}-\frac{\tilde g_j}{f_j}$ gives a \v{C}ech 1-cocycle with values in $\mathcal O_X$, 
 and its cohomology class is the image of $[\{g_i\}]$ by the map $H^0(\varphi(C), \mathcal N_{\iota})
 \to H^1(X, \mathcal O_X)$.
Note that its image in $H^1(X, \mathcal O_X^*)$ is represented by a \v{C}ech 1-cocycle 
 $\{\exp(\frac{\tilde g_i}{f_i}-\frac{\tilde g_j}{f_j})\}$.

Now take a divisor $D$ as in the claim. 
Let $\{h_i\}$ be the set of local defining functions of it on $\{U_i\}$.
Again, we set $h_i = 1$ when $U_i$ does not intersect $D$.
Then the meromorphic function $\frac{h_i}{f_i}$ gives a local frame of the sheaf $\mathcal O_X(\varphi(C)-D)$.
In particular, the functions $\frac{h_i}{f_i}(\frac{h_j}{f_j})^{-1}$ on $U_i\cap U_j$
 gives a \v{C}ech 1-cocycle with values in $\mathcal O_X^*$ which represents the class of 
 $\mathcal O_X(\varphi(C)-D)$.
If this class in contained in $\bar V$, then there is a set of local sections $\{g_i\}$ of $\mathcal O_X(\varphi(C))$
 such that the \v{C}ech 1-cocycle $\{\exp(\frac{\tilde g_i}{f_i}-\frac{\tilde g_j}{f_j})\}$
 gives the same class as $\{\frac{h_i}{f_i}(\frac{h_j}{f_j})^{-1}\}$.
It follows that there is a set of invertible holomorphic functions $\{\eta_i\}$ on $\{U_i\}$
 which satisfies
\[
\eta_i\eta_j^{-1} = \exp\left(\frac{\tilde g_i}{f_i}-\frac{\tilde g_j}{f_j}\right)
 \left(\frac{h_i}{f_i}\left(\frac{h_j}{f_j}\right)^{-1}\right)^{-1}
\]
 on $U_i\cap U_j$.
Therefore, we have 
\[
\exp\left(\frac{\tilde g_i}{f_i}-\frac{\tilde g_j}{f_j}\right)
 \left(\frac{\eta_ih_i}{f_i}\left(\frac{\eta_jh_j}{f_j}\right)^{-1}\right)^{-1} = 1
\]
 on $U_i\cap U_j$.
Since $\frac{\tilde g_i}{f_i}-\frac{\tilde g_j}{f_j} = 0$ on $U_i\cap U_j\cap \varphi(C)$, 
 we have 
\[
\frac{\eta_ih_i}{\eta_jh_j} = \frac{f_i}{f_j}
\]
 on $U_i\cap U_j\cap \varphi(C)$.
It follows that the set of functions $\{\eta_ih_i\}$ on $\{U_i\cap\varphi(C)\}$ gives a 
 section of $\mathcal N_{\iota}$.
On the other hand, since $\eta_i$ does not have zero, the set of functions $\{\eta_ih_i\}$ defines the divisor $D$.
It follows that the section $\{\eta_ih_i\}$ of $\mathcal N_{\iota}$ satisfies the property required in Theorem \ref{thm:CB}.\qed

\begin{rem}
Assume that the map $H^0(\varphi(C), \mathcal N_{\iota})\to H^1(X, \mathcal O_X)$ is surjective.
For example, this is the case when $\varphi(C)$ is sufficiently ample.
 then for each $p_i\in P$, if there is an effective divisor $D$ which is algebraically equivalent to $\varphi(C)$
 which avoids $p_i$ but passes through all points in $P\setminus\{p_i\}$, the map $\varphi$ is semiregular.
\end{rem}

\begin{example}
Consider a fibration $\pi\colon X\to B$ over a compact Riemann surface $B$
 which have a section $\sigma$.
Assume that $(n+1)\sigma(B)$ is very ample for some positive integer $n$.
Let $p_1, \dots, p_k$ be points on $B$ and $\sigma'$ be a smooth irreducible curve 
 in the linear equivalence class of $n\sigma(B) + \sum_{i=1}^k \pi^{-1}(p_i)$.
Assume that $\sigma(B)$ and $\sigma'$ intersect transversally
 and let $S = \sigma(B)\cap \sigma'$.
Let $S'=\{s_1, \dots, s_l\}$ be any subset of $S$, where $l\leq k+1$.
Let $C$ be the unchaining of $\sigma(B)\cup \sigma'$ obtained by normalizing 
 $S'$ and $\varphi\colon C\to X$ be the natural map.
Then $\varphi$ is semiregular, so unobstructed.
In fact, for any $j$, $1\leq j\leq l$, 
 take $D$ as $\sum_{i=1, i\neq j}^l\pi^{-1}(\pi(s_i)) + \bar\sigma$, 
 where $\bar\sigma$ is a general curve in the class $(n+1)\sigma(B)$.
This satisfies the condition of Corollary \ref{cor:geomCB}, since any two fibers are algebraically equivalent.
\end{example}

%Let $s_i$ be a section of $\mathcal O_X(\varphi(C))$ whose zero locus is the divisor in the statement.
%Then the restriction of it to $\varphi(C)$ gives a section of $\mathcal N_{\iota}$ whose
%associated first order deformation satisfies the condition required in Theorem \ref{thm:CB}.\qed

%\begin{example}
%Let $X$ be a smooth surface and 
% $C = C_1\coprod C_2$ be the disjoint union of nodal curves.
%Let $\psi\colon C\to X$ be a map such that the restriction of $\psi$ to $C_1$ and $C_2$ 
% is an embedding.
%Assume the image $\psi(C)$ is
% a reduced nodal curve and that 
% $\psi(C)$ is semiregular in the classical sense.
%
%Note that in this case, the embedded curve $C_1$ and $C_2$ are both semiregular in the classical sense
% if and only if the map $\psi$ is semiregular.
%Thus, we have a criterion for the semiregularity of these curves by　Theorem \ref{thm:CB}.
%Namely, let $P = \psi(C_1)\cap \psi(C_2)$ be the set of nodes.
%Then, the map $\psi$ is semiregular if and only if for each $p\in P$, there is a first order deformation of 
% $C$ which smoothes $p$ but does not smoothe the nodes in $P\setminus \{p\}$.

%\end{example}

\subsection{Criterion for semiregularity of reducible curves}
Let $X$ be a smooth projective surface.
Then $X$ contains many semiregular curves.
Namely, assume that a reduced curve $D\subset X$ satisfies the property that 
$H^1(X, \mathcal O_X(D)) = 0$.
In this case, the standard exact sequence
\[
0\to \mathcal O_X\to \mathcal O_X(D)\to \mathcal O_D(D)\to 0
\]
induces the cohomology exact sequence
\[
\cdots \to H^1(X, \mathcal O_X(D))\to H^1(D, \mathcal O_D(D))\to H^2(X, \mathcal O_X)\to\cdots.
\]
When  $H^1(X, \mathcal O_X(D)) = 0$,
the map $H^1(D, \mathcal O_D(D))\to H^2(X, \mathcal O_X)$ is injective.
Since this map is the dual of the semiregularity map $H^{0}(X, \mathcal K_X)\to 
H^{0}(D, \mathcal \iota^*K_X)$, 
it follows that $D$ is semiregular.
Note that the condition 
$H^1(X, \mathcal O_X(D)) = 0$ 
is satisfied when $D$ is sufficiently large.
For example, when $D$ is of the form $D = D' + K_X$, where $D'$ is nef and big, 
then by the Kawamata-Viehweg vanishing theorem \cite{Ka, V}, we have $H^i(X, \mathcal O_X(D)) = 0$ for $i = 1, 2$.

Corollary \ref{cor:CB} gives a criterion for unobstructedness of general nodal curves on surfaces.
When the curve has sufficiently high degree and multiple components, we can deduce a simpler criterion.
Let $C = \cup_{i=1}^m C_i$ be a connected reduced nodal curve,
 where $C_i$ is an irreducible component.
Let $\varphi\colon C\to X$ be an immersion.
%In the present case, the image need not be nodal nor non-reduced.
Assume that for each $1\leq i\leq n-1$, the intersection $C_i\cap C_{i+1}$ consists of 
 $\max\{K_X\cdot \varphi_i(C_i), K_X.\varphi_{i+1}(C_{i+1})\} + 1$ or more points.
 % (there can be other intersections
 %between $C_i$ and $C_j$, $i\neq j$).
Here $\varphi_i$ is the restriction of $\varphi$ to $C_i$.
In this situation, we have the following.
\begin{lem}\label{lem:union}
Assume that one of $\varphi_i$ is semiregular.
Then the map $\varphi$ is semiregular.
\end{lem}
\proof
We need to show that the map $H^0(X, \mathcal K_X)\to H^0(C, \varphi^*\mathcal K_X)$
is surjective.
By assumption, the map
$H^0(X, \mathcal K_X)\to H^0(C_i, \varphi_i^*\mathcal K_X)$
 is surjective.
The degree of the invertible sheaf $\varphi_i^*\mathcal K_X$ is $K_X.\varphi_i(C_i)$.
Then a section of $\varphi_i^*\mathcal K_X$ which is zero at $K_X.\varphi_i(C_i)+1$ points on $C_i$ must be zero.

Let $s$ be a section of $H^0(C, \varphi^*\mathcal K_X)$.
Its restriction to $C_j\cup C_{j+1}$ for any $j$
 gives an element of $H^0(C_j\cup C_{j+1}, \varphi^*\mathcal K_X|_{C_j\cup C_{j+1}})$.
Since $C_j\cap C_{j+1}$ consists of at least $K_X.\varphi_{j+1}(C_{j+1}) + 1$ points, 
 the value of the restriction $s|_{C_j}$ uniquely determines 
 $s_{C_{j+1}}$, and vice versa.
Therefore, for any $j$, the restriction $s|_{C_j}$ uniquely determines the entire section $s$.
Then, since the map $H^0(X, \mathcal K_X)\to H^0(C_i, \varphi_i^*\mathcal K_X)$ is surjective,
 the corresponding section $s$ must be also the pull back of a section of $H^0(X, \mathcal K_X)$. \qed\\

\subsection{Nodal curves with small genus on surfaces of non-negative Kodaira dimension} 
In this section, we prove the following theorem.
\begin{thm}\label{thm:2}
	Let $X$ be a smooth complex projective surface with an effective canonical class.
	Let $L$ be a very ample class.
	Then, there is a positive number $A$ which depends on $L$ such that 
	for any positive integer $n$, the numerical class of $nL$ contains
	an embedded irreducible nodal curve $C$ whose geometric genus is less than $An$.
\end{thm}
\proof
By replacing $L$ with $mL$ for some positive integer $m$, we may assume
%Let us take a class $L = K_X + D$, where $D$ is a very ample divisor.
%In particular, a smooth curve in the class $L$ is semiregular.
%By taking $D$ large enough, we assume 
 $L^2\geq K_X.L + 1$.
We also assume that $K_X.L+1\geq \dim H^0(X, \mathcal K_X)$.

Let $C_1, \dots, C_n$ be distinct curves in the class $L$ such that the union 
 $\cup_{i=1}^n C_i$ is an embedded nodal curve.
Let $C(l)$, $l = 1, \dots, n,$ be curves constructed as follows.
Namely, let $C(1) = C_1$.
Take $P$ to be a subset of $C_1\cap C_2$ which has cardinality $K_X. L+1$.
Then let $C(2)$ be the unchaining of $C_1\cup C_2$ obtained by normalizing
 the nodes except those in $P$.
Let $P'$ be a subset of $C_2\cap C_3$ which has cardinality $K_X.L+1$ and
 let $C(3)$ be the unchaining of $C_1\cup C_2\cup C_3$ obtained by normalizing
 the nodes except those in $P\cup P'$.
Repeating this, we construct a nodal curve $C(i)$ for all integers $i = 1, \dots, n$.
There is a natural map $\varphi(i)\colon C(i)\to X$ induced by the inclusions $C_i\hookrightarrow X$.
By Lemma \ref{lem:union}, the map $\varphi(i)$ is semiregular for any $i$.

Now we prove the following.
\begin{lem}
There is a first order deformation of $\varphi(n)$ which simultaneously smoothes some of the nodes
 $C_i\cap C_{i+1}$ for each $i = 1, \dots, n-1$.
\end{lem}
\proof
First consider the map $\varphi(2)\colon C(2)\to X$.
Take the exact sequence of sheaves 
\[
0\to \mathcal N_{\varphi(2)}(-P)\to \mathcal N_{\varphi(2)}\to \oplus_{p\in P}\Bbb C_p \to 0
\]
 on $C(2)$,
 where $\mathcal N_{\varphi(2)}$ is the normal sheaf of the map $\varphi(2)$ and 
 $\Bbb C_p$ is the skyscraper sheaf supported at $p$.
The associated cohomology exact sequence is
\[
\begin{array}{ll}
0\to H^0(C(2), \mathcal N_{\varphi(2)}(-P))  \to H^0(C(2), \mathcal N_{\varphi(2)})
 & \to \oplus_{p\in P}\Bbb C_p\\
 & \to
H^1(C(2), \mathcal N_{\varphi(2)}(-P))\to H^1(C(2), \mathcal N_{\varphi(2)})\to 0.
\end{array}
\]
Note that the sheaf $\mathcal N_{\varphi(2)}(-P)$ is isomorphic to $\pi_*\mathcal N_{\psi}$, 
where $\pi\colon C_1\coprod C_2\to C(2)$ and
$\psi\colon C_1\coprod C_2\to X$ are natural maps.
Then by the Leray spectral sequence, we have 
\[
H^i(C(2), \mathcal N_{\varphi(2)}(-P))\cong 
H^i(C_1\coprod C_2, \mathcal N_{\psi})
\cong 
H^i(C_1, \mathcal N_{C_1})\oplus H^i(C_2, \mathcal N_{C_2}), \;\; i = 0, 1,
\] 
where $\mathcal N_{C_1}$ and $\mathcal N_{C_2}$ are the normal sheaves of $C_1$ and $C_2$, respectively.
In particular, the dual of $H^1(C(2), \mathcal N_{\varphi(2)}(-P))$ is isomorphic to 
\[
(H^1(C_1, \mathcal N_{C_1})\oplus H^1(C_2, \mathcal N_{C_2}))^{\vee}\cong 
H^0(C_1, \mathcal K_X|_{C_1})\oplus H^0(C_2, \mathcal K_X|_{C_2})
\]
 by the Serre duality.
On the other hand, the dual of $H^1(C(2), \mathcal N_{\varphi(2)})$
is isomorphic to $H^0(C(2), \varphi(2)^*\mathcal K_X)$.
Therefore, the dual of the last three terms of the above exact sequence is
equivalent to 
\[
0\to H^0(C(2), \varphi(2)^*\mathcal K_X)\to H^0(C_1, \mathcal K_X|_{C_1})\oplus H^0(C_2, \mathcal K_X|_{C_2})
\to \oplus_{p\in P}\Bbb C_p.
\]

Recall that $C_1, C_2$ are semiregular in the classical sense, and the map $\varphi(2)$ is also semiregular.
In particular, we have $H^0(C(2), \varphi(2)^*\mathcal K_X)\cong H^0(X, \mathcal K_X)$
 and $H^1(C_1, \mathcal K_X|_{C_1})\oplus H^0(C_2, \mathcal K_X|_{C_2})\cong H^0(X, \mathcal K_X)^{\oplus 2}$.
Since we assumed that the cardinality of $P$ is larger than $\dim H^0(X, \mathcal K_X)$, 
 it follows that the map $H^0(C_1, \mathcal K_X|_{C_1})\oplus H^0(C_2, \mathcal K_X|_{C_2})
 \to \oplus_{p\in P}\Bbb C_p$
 is not surjective.
In the original exact sequence, this means that the map $\oplus_{p\in P}\Bbb C_p\to 
 H^1(C(2), \mathcal N_{\varphi(2)}(-P))$ is not injective, which in turn implies that
 the map $H^0(C(2), \mathcal N_{\varphi(2)})
 \to \oplus_{p\in P}\Bbb C_p$ is not zero.
Therefore, there is a section $s$ of the normal sheaf $\mathcal N_{\varphi(2)}$
 which corresponds to a first order deformation of $\varphi(2)$ which smoothes some of the nodes in $P$.

Now consider the union $C(2)' = C_2\cup C_3$ and the map $\varphi(2)'\colon C(2)'\to X$
 obtained by restricting $\varphi(n)$.
By the same argument as above, there is a section $s'$ of the normal sheaf $\mathcal N_{\varphi(2)'}$
 of $\varphi(2)'$ which corresponds to 
 a first order deformation of $\varphi(2)'$
 which smoothes some of the nodes in $C_2\cap C_3$.
On the other hand, the restriction of $\mathcal N_{\varphi(2)'}$ to $C_2$ is naturally isomorphic to the 
 restriction of the sheaf $\mathcal N_{\varphi(2)}(-P+P')$ to $C_2$.
Similarly, the restriction of $\mathcal N_{\varphi(2)}$ to $C_2$ is naturally isomorphic to 
 the 
 restriction of the sheaf $\mathcal N_{\varphi(2)'}(P-P')$ to $C_2$.
In particular, we can add $s$ and $s'$ to obtain a section of the restriction of the normal sheaf 
 $\mathcal N_{\varphi(n)}$ to $C(3)$.
Repeating this argument, we obtain a desired first order deformation.\qed\\

By Theorem \ref{thm:main}, such a deformation can be extended to an algebraic deformation 
 $\tilde\varphi(n)\colon \tilde C(n)\to X$, where $\tilde C(n)$ is a deformation of $C(n)$.
Since the map $\varphi(n)$ is an immersion and the image is a nodal curve, 
 we can also assume the same property to $\tilde\varphi(n)$.
The geometric genus of the curve $\tilde C(n)$ is at most $ng(C_1) + (n-1)K_X.L$.
Therefore, the claim follows.\qed

\begin{rem}
	\begin{enumerate}
		\item In many cases (e.g., when $L$ is sufficiently large),
		 the constant $A$ in the theorem can be taken as $A = g(C_1) + K_X.L$.
		 
		\item Note that a general smooth curve in the numerical class of $nL$ has genus of order
		$\frac{L^2}{2}n^2$.
		 %Therefore, the curve obtained in the theorem has very large number of nodes.
		In particular, by taking $A$ slightly bigger if necessary, 
				 we can assume that the number of nodes $\delta(C)$ of $C$ satisfies
				 $A\delta(C) > g(C)^2$, where $g(C)$ is the geometric genus of $C$.
		\item According to the proof, the surface $X$ is dominated by a family of nodal curves 
		which satisfy the property in the theorem.
	\end{enumerate}
\end{rem}

\section*{Acknowledgement}
\noindent
The author was supported by JSPS KAKENHI Grant Number 18K03313.

\end{document}